\documentstyle[leqno,amssymb,
amstex,12pt]{amsart}
\normalsize
\newcommand {\e}		{\varepsilon}

\numberwithin{equation}{section}

\newtheorem{th}		          	{Theorem}

\newtheorem{lemma}[equation]		{Lemma}
\newtheorem{prop}[equation]		{Proposition}
\newtheorem{corollary}			{Corollary}

\theoremstyle{definition}
\newtheorem{definition}[equation]	{Definition}

\newtheorem{examples}[equation]		{Examples}

\theoremstyle{remark}
\newtheorem{remark}[equation]		{Remark}

\begin{document}

\centerline{\bf GROWTH AND SPECTRUM OF DIFFUSIONS}

\vspace{18pt}
\centerline{\bf Lino Notarantonio\footnote{The author is a Dov
Biegun Postdoctoral Fellow at the Weizmann Institute of Science}}
\centerline{Department of Theoretical Mathematics}
\centerline{The Weizmann Institute of Science,} 
\centerline{Rehovot, 76100, Israel}



\begin{abstract}
We prove an upper bound on the bottom of the essential spectrum of a
diffusion in term of the growth of the volume of $X$, generalizing a
result by R.~Brooks \cite{brooks1}.
\end{abstract}


\thispagestyle{empty}

\section*{Introduction}
Let $X$ be a locally compact, second countable Hausdorff
space $X$ and let $dx$ be a Radon measure on $X$, which
sometimes we shall call the volume of $X$.

Let \( H\ge 0\) be a self-adjoint operator 
associated with a diffusion $(a, D[a])$, that is, a 
strongly local, regular Dirichlet form on $L^2 = L^2(\Omega,dx)$
\cite{fukufuku}. In this paper we examine 
the relation between the growth (of the volume) of $X$ and the bottom
of the spectrum of the operator $H$. More
precisely, we introduce a (pseudo-)distance $d(x,y)$ on $X$ through
$(a,D[a])$ and assume that the metric topology is equivalent to the
original topology of $X$; we fix a point $o \in X$ and, denoting by 
$V(r)$ the measure of the ball $B(o,r)$, $r>0$, introduce the
growth of $X$ by means of the quantity
\[
\mu := \limsup_{r \to +\infty} \frac{\ln V(r)}{r};
\]
finally we also introduce 
\begin{equation}
\lambda_e := \sup \lambda_o(X \setminus K),
\label{eqn:bottom-ess}
\end{equation}
where the supremum is taken on the family of compact subsets of $X$,
and 
\[ 
\lambda_o(X\setminus K) = \inf\left\{ 
\frac{a[u,u]}{\int_X u^2 dx} : u\not\equiv 0,\ u \in D[a]\cap C_c(X
\setminus K)
\right\}
\]
is the bottom of the spectrum of \( H\) on $L^2(X\setminus K)$ with
Dirichlet boundary condition on $\partial K$. 

The formula (\ref{eqn:bottom-ess}) above can be seen as a
generalization of the well-known formula \cite{donnelly}, 
\cite{donnelly-li} that holds when $H$ is the
Laplace-Beltrami operator on a Riemannian manifold $X$ ($dx$ is then
the Riemannian volume) and $\lambda_e$ is the bottom of the essential
spectrum of $H$. The number $\lambda_e$ may be infinite, which is the
case when $H$ has discrete spectrum and the essential spectrum is
empty. 

The main result of this paper is an extension of a result by
R.~Brooks \cite[Theorem~1]{brooks1} and goes as follows. 

\begin{th}\label{thm:main-result}
If $X$ is non-compact, the volume of $X$ is infinite and the metric
space $(X,d)$ is complete, then 
\[
\lambda_e \le \frac{\mu^2}{4}.
\]
\end{th}

In particular, if the measure $dx$ has polynomial growth (e.g., if
$dx$ satisfies a doubling condition), then $\mu=0$, hence the bottom
of the essential spectrum is $\lambda_e=0$. See
Corollary~\ref{corollary:main-result}. 
On the other hand, if the growth is exponential (as in the case, e.g.,
of $H^{n+1}$, the $n+1$-dimensional hyperbolic space), then  we
recover a well-known result on the bottom of the Laplace-Beltrami
operator on $H^{n+1}$; cf., e.g., the result by H.P.~McKean
\cite{mckean}.

\vspace{.0625in}
The above result sharpens a similar one for \( \lambda_o(X)\) (with \(
K = \emptyset\)) by K.-T.~Sturm \cite[Theorem~5]{sturm1}.

\vspace{.0625in}
Theorem~\ref{thm:main-result} is sharp in the following sense. 
It is known that non-complete metric spaces may have discrete spectra;
cf. e.g. results for the Laplace operator in ``quasibounded
domains'' of $R^n$ by D.~Hewgill \cite{hewgill1}, \cite{hewgill2} (see
also the references cited in these two papers). Moreover, 
compact Riemannian manifolds and some Riemannian manifolds with finite
volume are known to have discrete spectrum: the former is classical,
while instances of the latter have been constructed by H.~Donnelly \&
P.~Li in \cite{donnelly-li}.

\vspace{.0625in}
On the other hand, it must be remarked that there are examples where
the upper bound on 
$\lambda_e$ in term of the growth of $X$ is not sharp in that
$\lambda_e=0$ and $\mu>0$; cf. the discussion in the
Introduction of \cite{brooks1}.

\vspace{.0625in}
The organization of the paper is as follows: in the first section we
fix the notation, introduce the relevant concepts and prove the
preliminary results that will be needed in the second section, which is
devoted to the proof of the main result. The proof of our main result
follows the lines of the proof of Theorem~1 in \cite{brooks1} and is
based on Theorem~\ref{thm:agmon}, which is a result that may have some
interest in its own.

\section{Preliminaries}\label{sec:preliminaries}

Let $X$ be a locally compact, second countable Hausdorff space
and let $dx$ be a Radon measure on it. 

\vspace{.0625in}
{\em Dirichlet forms~\cite{fukufuku}.}
We let $(a, D[a])$ denote the (Dirichlet) form associated with the
self-adjoint operator $H$ on $L^2 = L^2(X,dx)$, so that
\[
\langle Hf, g\rangle = a[f,g],\ \ \ f\in D(H), \ \ g\in D[a], 
\]
where \( \langle \cdot,\cdot\rangle\) denotes the inner product in \(
L^2\).  

We shall say that \( u\in D_{loc}[a]\) if for every compact set $C
\subset X$ there exists $\overline u \in D[a]$ such that $u =
\overline u$, $dx$-a.e. on $C$.

We shall assume in the sequel that $(a,D[a])$ is a diffusion, {\em
i.e.,\/} $(a,D[a])$ satisfies the following property: 
\[
a[u,v]=0,\ \ \ u,v\in D[a], 
\]
whenever $u=\mbox{const.}$ on $\mbox{supp}\,v$. 

Then it is standard that we can write the form as 
\[
a[u,u] = \int_\Omega d\Gamma[u,u],
\]
where the map $(u,v) \mapsto d\Gamma[u, v]$, defined on $D[a] \times
D[a]$ with values in the space of signed Radon measures, is a
non-negative symmetric bilinear form (the energy measure associated
with the form $(a,D[a])$). The energy measure can be defined
according to the formula 
\begin{equation}
\int_X \phi\, d\Gamma[u,u] = a[u,\phi u] - \frac{1}{2} a[u^2,\phi],
\label{eqn:energy-measure}
\end{equation}
for every $u \in D[a] \cap L^\infty(X,dx)$ and every $\phi \in D[a]
\cap C_c(X)$.

\begin{prop}\label{prop:energy-measure}
The energy measure satisfies the following properties:

\begin{itemize}

\item[(D$_1$)] {\em The Leibnitz rule.\/} For every $u,v \in
D_{loc}[a] \cap L^\infty(X,dx)$ and $w \in D_{loc}[a]$ 
\[
d\Gamma [uv,w] = v(x ) d\Gamma [u,w] + 
u(x)d\Gamma [ v , w]
\]
in the sense of measures.

\item[(D$_2$)] {\em The Schwarz rule.\/} 
If $u,v \in D[a]$, $f \in L^2(\Omega, d\Gamma [u,u])$, $g\in L^2
(\Omega , d\Gamma [v,v])$, then $fg$ is integrable w.r.t. the
absolute variation $|d\Gamma [u,v]|$ of $d\Gamma [u,v]$, and
\[
|fg| \, |d\Gamma [u,v]| \le 
\frac{\eta}{2}|f|^2 d\Gamma [u,u] + \frac{1}{2\eta} |g|^2
d\Gamma [v,v],
\]
for $\eta >0$;

\item[(D$_3$)]  {\em The chain rule.\/} Let $\eta \in C^1(R)$ with
bounded derivative. Then $u \in D_{loc} [a]$ implies $\eta(u) \in D[a]$
and \[
d\Gamma [\eta(u), v] = \eta'(u) d\Gamma[u,v],
\]
for every $v\in D[a] \cap L^\infty(X,dx)$. 

\item[(D$_4$)] {\em The truncation property.\/}
(\cite[Lemma~3-(o)]{mosco1} For every $u\in D_{loc}[a]$
\[
d\Gamma[u_+ , v] = 1_{\{ u >0\}}  d\Gamma[u,v],
\]
where $u_+(x) := \max\{ u(x), 0\}$, $x \in X$.

\item[(D$_5$)] {\em Locality.\/} If $A$ is an open set and $u_1 = u_2$
$dx$-a.e. on $A$, $u_1,u_2\in D[a]$, then 
\[
1_A(x) d\Gamma[u_1,u_1] = 1_A(x) d\Gamma[u_2,u_2]
\]
on $X$; moreover 
\[
1_A(x) d\Gamma[u,v] = 0,
\]
on $X$, whenever $u\in D[a]$ is constant on $A$, for arbitrary $v \in
D[a]$.  
\end{itemize}
\end{prop}

\begin{definition}\label{ccd}
For $x, y \in \Omega$ let 
\[
d(x,y) := \sup \{ \psi (x) - \psi (y) : \psi \in C,\ \ d\Gamma [\psi,
\psi] \le dx\}.
\]
Then it is not difficult to prove that $d(\cdot , \cdot)$ is a
(pseudo-)distance on $\Omega$, which we call the Carnot-Carath\'eodory
distance (associated with the form $(a, D[a])$);
cf. \cite[\S~III-4]{vscc}. 
\end{definition}

\begin{remark}\label{remark:ccd}
Notice that $d(x,y)=0$ not necessarily
implies $x=y$. Moreover $d(x,y)$ may be equal to $0$ or $\infty$,
for some  $x\not= y$.
\end{remark}

\vspace{.25in}
{\em Warning.\/}~In the rest of the paper we shall make the assumption
that the 
topology induced by the metric is equivalent to the original topology
of $X$ and that $(X,d(\cdot,\cdot))$ is a complete metric space.

This in turn is equivalent \cite{sturm6} (cf. also \cite[\S
4]{sturm1}) to the fact that all balls are relatively compact in $X$.

\begin{remark}\label{reamrk:di}
Notice that in the case of the Dirichlet integral on a bounded open
set of $R^n$ the energy measure $d\Gamma[u,u]$ is absolutely
continuous with respect to the Lebesgue measure; the Radon-Nikodym
derivative is equal to $|\nabla u|^2$ so that
\[
d(x,y) = \sup \{ \psi (x) - \psi (y) : \psi \in C,
\ \ |\nabla \psi| \le 1\ \mbox{a.e.} \}.
\]
Notice that $d(x,y)$ is locally equivalent to the standard euclidean one,
that is, for any $x\in \Omega$ and any neighborhood $U$ of $x$, there
exists a constant $c >0$ such that  
\[
c^{-1}\, d(x,y) \le |x-y| \le c\,d(x,y),
\]
for $y\in U$.
\end{remark}

We shall need the following technical result to localize functions in
$D[a]$.

\begin{lemma}\label{lemma:distance}
For each compact set $K \subset X$ there exists a function $\chi$ such
that $\chi_K \in D_{loc} \cap C(X)$, $\chi_K(x) =1$ on $K$, the
support of $\chi_K$ is contained in a neighborhood $B$ of $K$ and
$d\Gamma[\chi_K,\chi_K] \le 16 (\mbox{diam}\, B)^{-2} dx$.  
\end{lemma}

\begin{pf} As $K$ is compact, then $K \subset B(x_o,R)$, for some
$x_o\in X$, $R>\,\mbox{diam}\,K \ge 0$, and $B = B(x_o0,2R)$. 
Consider the function
$\eta \in C^1(R)$ such that $\eta(t) =1$, for $t \in (-\infty,1)$,
$\eta(t)=0$, for $t \in [2,+\infty)$ and $|\eta'| \le 1$; let 
\[
\chi_K(x) := \eta\left( \frac{d(x,x_o)}{R}\right).
\]
The function $\rho_{x_o}(x):= d(x,x_o)$ is in $D_{loc}[a] \cap C(X)$
and $d\Gamma[ \rho_{x_o}, \rho_{x_o}] \le dx$ 
\cite[\S~4, Lemma~A$^\prime$]{sturm1}, thus by the chain rule $\chi_K \in
D_{loc}\cap C(X)$ and $d\Gamma[\chi_K,\chi_K] \le R^{-2} dx$;
moreover, by definition, $\chi_K =1$ on ($B(x_o,R)$, hence on) $K$ and
$\chi_K = 0$ outside $B = B(x_o,2R)$, so the proof is
completed. \end{pf}  

\begin{prop}\label{prop:distance}
Let $A$ be any subset of $X$; then 
\[
\phi(x) := \mbox{dist}\,(x,A) = \inf\{ d(x,y) : y\in A\}
\]
is in $D_{loc}[a]\cap C(X)$, with $d\Gamma[\phi,\phi] \le dx$. 
\end{prop}

\begin{pf}
We first prove the continuity of the map $\phi$, by showing that 
$\phi(x_n) \rightarrow \phi(x)$ whenever $d(x,x_n) \rightarrow 0$. By
definition \[
\phi(x_n) \le d(x_n,y),\ y\in A,
\]
thus 
\[
\limsup_{n\to +\infty} \phi(x_n) \le \limsup_{n\to +\infty}
d(x_n,y) =d(x,y),
\]
hence, taking the supremum over all $y \in A$, 
\[
\limsup_{n\to+\infty} \phi(x_n) \le \inf\{ d(x,y): y \in A\} =
\phi(x). 
\]
On the other hand for any $\e>0$ there is $y_\e \in A$ such that 
\[
d(x_n,y_\e) \le \phi(x_n) + \e,
\]
thus. choosing $\e = 1/n$ and using the triangle inequality
$d(x,y_\e) \le d(x,x_n) + d(x_n,y_\e)$, we have 
\[
\phi(x) \le  d(x,x_n) + \phi(x_n) + 1/n;
\]
hence
\begin{align*}
\phi(x) & \le  \liminf_{n\to +\infty} \Big[ d(x,x_n) + \phi(x_n) + 
1/n \Big] \\
& = \liminf_{n\to +\infty} \Big[ \big( d(x,x_n) + 1/n\big) + 
\phi(x_n) \Big] \\
& \le \limsup_{n \to +\infty} \Big(d(x_n,x) + 1/n \Big) + 
\liminf_{n\to +\infty} \phi(x_n) \\
& \le \liminf_{n\to +\infty} \phi(x_n)
\end{align*}
thus 
\[
\phi(x) \le  \liminf_{n\to +\infty} \phi(x_n).
\]
Therefore 
\[
\limsup_{n\to+\infty} \phi(x_n) \le 
\phi(x) \le \liminf_{n\to +\infty} \phi(x_n),
\]
which implies that the map $\phi$ is continuous on $X$. 

Now let us prove that $\phi \in D_{loc}[a]$, with
$d\Gamma[\phi,\phi] \le dx$. Again by definition there exists a
sequence $(y_n)$ of points in $A$ such that 
\[
\phi(x) = \lim_{n \to +\infty} d(x,y_n).
\]
Consider the map $\phi_n(x) := d(x,y_n)$, $n = 1,2,\ldots$, so that
$\phi_n(x) \rightarrow \phi(x)$, $x\in X$; in fact, the convergence is
uniform on any relatively compact open subset $Y \subset X$. Indeed,
it is easy  to see that each $\phi_n$ is Lipschitz continuous with
Lipschitz constant less than or equal to 1 so that the sequence is
pre-compact in $C(Y)$, according to the Ascoli-Arzel\'a
criterion. Moreover and without loss of generality, we can also assume
that $\phi_n \le \phi_{n+1}$, $n=1,2,\ldots$ (consider otherwise
$\sup_{1\le i\le n} \phi_n$). We notice  
\cite[\S~4, Lemma~A$^\prime$]{sturm1} that 
$\phi_n \in D_{loc}[a] \cap C(X)$ and $d\Gamma[\phi_n,\phi_n] \le dx$,
$n=1,2\ldots$. By localization on $Y$, we also have that  
\[
\int_Y d\Gamma[\phi_n,\phi_n] + \int_Y \phi_n^2 dx \le 
\mbox{const}(Y).
\]
Therefore the family 
\[
\{ \phi_n : \phi_n \in D_{loc}[a] \cap C(X),\ \ d\Gamma[\phi_n,\phi_n]
\le dx,\ \ n\in N\}
\]
is convex and uniformly bounded in $D[a]$. By the Banach-Saks
theorem, there exists (possibly a subsequence of)
$(\phi_n)$) that converges strongly (in $D[a]$, hence strongly in)
$L^2$ to some $\overline \phi \in D_{loc}[a]$; from the strong
convergence in $D[a]$ we also have $d\Gamma[\overline \phi,\overline
\phi] \le dx$ as measures.
On the other hand the whole sequence converge uniformly on $Y$ to
$\phi$, therefore $\phi(x) = \overline \phi (x)$ for q.e. $x\in
Y$, and so we have $\phi \in D_{loc}[a]$ and also $d\Gamma[\phi,\phi]
\le dx$.
\end{pf}

\section{The main result}

Denote by \( \rho(x) := d(x,o)\), the function ``distance from a given
point $o\in X$''.

As in Brooks's paper \cite{brooks1}, our Theorem~\ref{thm:main-result}
is a consequence of the following result. 

\begin{th}\label{thm:agmon}
Let $K \subset X$ be a compact set (possibly empty), let
\[
\lambda_o(X \setminus K) = \inf
\left\{ \frac{a[u,u]}{\displaystyle{\int_{X \setminus K} u^2 dx}} : 
0\not = u\in D[a] \cap C_c(X \setminus K)
\right\}.
\]
If 
\[
\int_{X \setminus K} \exp(-2\alpha \rho(x))\,dx <+\infty, 
\]
for some \( \alpha \in (0, \sqrt{\lambda_o(X \setminus K)\,})\), then 
\[
\int_{X \setminus K} \exp(2\alpha \rho(x))\,dx <+\infty.
\]
\end{th}

\begin{pf*}{Theorem~\ref{thm:agmon} implies
Theorem~\ref{thm:main-result}} 
Recall that $V(r)$ stands for the volume of $B(o,r)$, $r>0$. 
If $2\alpha < \mu$ then 
\begin{align*}
\int_{X \setminus K} \exp(-2\alpha \rho(x))\,dx & \le
\sum_{r =1}^{+\infty} \Big[ V(r) - V(r-1)\Big]e^{-2\alpha(r-1)} \\
& = \sum_{r =1}^{+\infty} V(r)e^{-2\alpha r}\Big[ e^{2\alpha} -1\Big]
\end{align*}
and the latter sum is finite, as it follows by comparing it with a
geometric series and by the fact that $2\alpha > \mu$. Thus if
$2\alpha > \mu$ and $\alpha < \sqrt{\lambda_o(X \setminus K)\,}$ then
by Theorem~\ref{thm:agmon} it follows that 
\[
\int_{X \setminus K} \exp(2 \alpha\rho(x))\,dx < +\infty,
\]
but this inequality cannot be true, as the volume of ($X$, hence the
volume of) $X \setminus K$ is infinite. (Recall that $dx$ is a Radon
measure, hence the volume of any compact set is finite.)
Therefore $2\alpha \le \mu$, and letting \( \alpha\) approach
$\sqrt{\lambda_o(X \setminus K)\,}$, we have 
\[
\lambda_o(X \setminus K) \le \frac{\mu^2}{4},
\]
and the right-hand side does not depend on $K$. Taking the supremum
over $K$, we have \( \lambda_e \le \mu^2/4\), 
which proves Theorem~\ref{thm:main-result}.
\end{pf*}

\begin{pf*}{Proof of Theorem~\ref{thm:agmon}} Without loss of
generality we can assume that $\lambda_o( X\setminus K) \not=0$ and
hence, by a possible rescaling of the measure $dx$, that
$\lambda_o(M\setminus K)=1$. Let as consider the function $f(x) =
e^{h(x)}\chi(x)$, where $\chi \in D[a] \cap C_c(X \setminus
K)$ and $h$ is a bounded function in $D[a]$. (We'll make a
choice of these two functions later on.) 
Let us compute $d\Gamma[f,f] = d\Gamma[\chi e^h, \chi e^h]$. We have 
\begin{align*}
d\Gamma[f,f] & = d\Gamma[\chi e^h, \chi e^h] \\
\mbox{(Leibniz rule)}\ \ \ \ \  
& = e^{2h}d\Gamma[\chi,\chi] + \chi^2 d\Gamma[e^h,
e^h] + 2\chi e^h d \Gamma[e^h,\chi] \\
\mbox{(chain rule)}\ \ \ \ \  
& = e^{2h}d\Gamma[\chi,\chi] + e^{2h}\chi^2
d\Gamma[h, h] + 2\chi e^{2h} d \Gamma[h,\chi] \\
& = e^{2h}d\Gamma[\chi,\chi] + f^2 d\Gamma[h, h] + 
2\chi e^{2h} d \Gamma[h,\chi].
\end{align*}
As 
\[
1 = \lambda_o(X \setminus K) \le \frac{\displaystyle{\int_{X \setminus
K} d\Gamma[f,f] }}{\displaystyle{
\int_{X \setminus K} f^2\,dx}},
\]
we have
\[
\int_{X \setminus K} f^2\,dx \le \int_{X \setminus K}
\Big[ e^{2h} d\Gamma[\chi,\chi] + f^2 d\Gamma[h,h] + 2\chi e^{2h}
d\Gamma[h,\chi] \Big], 
\]
that is, 
\begin{align}
\int_{X \setminus K} f^2\,dx  & - \int_{X \setminus K} f^2 d\Gamma[h,h]
\label{eqn:pf-main-result1} \\
& \le \int_{X \setminus K} \Big[ e^{2h} d\Gamma[\chi,\chi] + 2\chi e^{2h}
d\Gamma[h,\chi] \Big].  
\notag
\end{align}
Moreover, by the Schwarz rule and taking into account that \( \chi \)
is a bounded function, 
\begin{align}
\int_{X \setminus K} 2\chi e^{2h} d\Gamma[h,\chi] & \le
C \left( \int_{X \setminus K} e^{2h} d\Gamma[\chi,\chi] \right)^{1/2}
\label{eqn:pf-main-result2} \\
&\ \ \ \times \left( \int_{X \setminus K} e^{2h} d\Gamma[h,h]
\right)^{1/2}, \notag
\end{align}
for some constant $C>0$. Now we turn to the choice of the functions
$h$, $\chi$. 
Let $(K_i)$ be a sequence of compact sets in $X \setminus
K$ which increases  to $X \setminus K$. For $i=1,2,\ldots$, let
\[
\chi_i(x) := \begin{cases}
\displaystyle{
\frac{1}{\delta} \mbox{dist}\,(x, X \setminus K_i)}, & 
\mbox{if $0 \le \mbox{dist}\,(x, X \setminus K_i) \le \delta$} \\
\\
1, & \mbox{if $ \mbox{dist}\,(x, X \setminus K_i) > \delta$.} 
\end{cases}
\]
By Proposition~\ref{prop:distance} $\chi_i \in D_{loc}\cap C(X)$,
$d\Gamma [\chi_i,\chi_i] \le \delta^{-2} dx$ and
by definition the support of $\chi_i$ is compact and contained in $X
\setminus K$; moreover, by the locality property of the energy
measure, the support of $d\Gamma[\chi_i,\chi_i]$ is contained in a
neighborhood $B_\delta(\partial K_i)$ of $K_i$. Furthermore, let us
choose the function $h$
such that $d\Gamma[h,h] \le \alpha\,dx$, for $\alpha \in
(0,1)$. Then from (\ref{eqn:pf-main-result1}),
(\ref{eqn:pf-main-result2}) we get  
\[ 
\int_{X \setminus K} f^2 (1-\alpha^2) dx \le C\left(\frac{2}{\delta} +
\frac{1}{\delta^2}\right) \int_{ B_\delta(\partial K_i)} e^{2h} dx. 
\]
Now the proof follows the proof given by R.~Brooks in
\cite[Theorem~2]{brooks1}: consider 
\[
h_n(x) := \min\{ \alpha \rho(x), -\alpha\rho(x) + n \}, \ n=1,2\ldots;
\]
notice that $h_n \in D_{loc}[a] \cap L^\infty(X,dx)$, $n=1,2\ldots$; 
$d\Gamma [h_n,h_n] \le \alpha\,dx$ and, under the assumption that
$\exp(-\alpha \rho(x))$ is integrable on $X \setminus K$, $h_n$ is
integrable for all $n$, so $h_n$ is an admissible function in the
definition of $f$. Notice also that $h_n(x) \le h_{n+1}(x) \rightarrow
h=\alpha\, \rho(x)$, $x \in X\setminus K$, as $n \to +\infty$. Thus
for $n$ sufficiently large,  
\[
\int_{X\setminus K} e^{2h_n}(1 -\alpha^2)\,dx \le
C\left(\frac{2}{\delta} + 
\frac{1}{\delta^2}\right) \int_{ B_\delta(\partial K)} e^{2\alpha
\rho(x)} dx, 
\]
that is, 
\[
\int_{X\setminus K} e^{2h_n} dx \le \mbox{const}
<+\infty, 
\]
where the constant at the right-hand side of the above inequality does
not depend on $n$. Taking the limit, as $n \to +\infty$, we have 
\[
\int_{X\setminus K} \exp(2\alpha \rho(x)) dx \le \mbox{const}
<+\infty, 
\]
which concludes the proof of the theorem.
\end{pf*}

\begin{corollary}\label{corollary:main-result}
If $X$ has sub-exponential growth, then $\lambda_e =0$. 
\end{corollary}
The condition that $X$ has sub-exponential growth is precisely the
fact that $\mu=0$.

\vspace{.25in}
\centerline{ACKNOWLEDGMENTS.}
It is a pleasure to thank M.~Solomyak for many
discussions and suggestions during the preparation of this paper; also
a discussion on a preliminary version of the paper with G.V.~Rosenblum
(who, together with M.~Solomyak, brought to the author's attention the
papers by D.~Hewgill) is gratefully acknowledged.

\bibliographystyle{plain}
\bibliography{/home/notaran/tex/references}

\end{document}